\documentclass{amsart}

\usepackage[colorlinks,driverfallback=dvipdfm]{hyperref}
\usepackage{color,graphicx,shortvrb}
\usepackage[utf8]{inputenc}

\usepackage[active]{srcltx} %SRC Specials for DVI search

\usepackage{enumerate}
\usepackage{amssymb}
\usepackage{tikz-cd}

\usepackage{parskip}

\usepackage{faktor}

\usepackage{ bbm }

\usepackage{comment}

\usepackage{amsmath,amssymb,verbatim,amscd,bm,color,yhmath,bm, enumerate}
\usepackage[normalem]{ulem}
\usepackage[driverfallback=dvipdfm]{hyperref}

\usepackage[active]{srcltx} %SRC Specials for DVI search

\newtheorem{lem}{Lemma}[section]
\newtheorem{thm}[lem]{Theorem}
\newtheorem{prop}[lem]{Proposition}
\newtheorem{cor}[lem]{Corollary}

\theoremstyle{definition}

\newtheorem*{ex}{Example}
\newtheorem{rem}[lem]{Remark}

\numberwithin{equation}{section}

 \DeclareSymbolFont{largesymbol}{OMX}{yhex}{m}{n}
 \DeclareMathAccent{\Widehat}{\mathord}{largesymbol}{"62}
 \DeclareSymbolFont{largesymbol}{OMX}{yhex}{m}{n}
 \DeclareMathOperator{\supp}{supp}
 \DeclareMathOperator{\coz}{coz}

\usepackage{ulem}
\usepackage{marginnote}
\usepackage{geometry}

\begin{document}

%\baselineskip24pt

%%%%% Title %%%%%

\title[Linear orthogonality preservers between commutative JB$^*$-triples]
{Linear orthogonality preservers between function spaces associated with commutative JB$^*$-triples}

\author[D. Cabezas]{David Cabezas}
\address[D. Cabezas]{
Departamento de An{\'a}lisis Matem{\'a}tico, Facultad de
Ciencias, Universidad de Granada, 18071 Granada, Spain.}
\email{dcabezas@ugr.es}

\author[A.M. Peralta]{Antonio M. Peralta}
\address[A.M. Peralta]{Instituto de Matem{\'a}ticas de la Universidad de Granada (IMAG), Departamento de An{\'a}lisis Matem{\'a}tico, Facultad de
Ciencias, Universidad de Granada, 18071 Granada, Spain.}
\email{aperalta@ugr.es}


\begin{abstract} It is known, by Gelfand theory, that every commutative JB$^*$-triple admits a representation as a space of continuous functions of the form $$C_0^{\mathbb{T}}(L) = \{ a\in C_0(L) : a(\lambda t ) = \lambda a(t), \ \forall \lambda\in \mathbb{T}, t\in L\},$$ where $L$ is a principal $\mathbb{T}$-bundle and $\mathbb{T}$ denotes the unit circle in $\mathbb{C}.$ We provide a description of all orthogonality preserving (non-necessarily continuous) linear maps between commutative JB$^*$-triples. We show that each linear orthogonality preserver  $T: C_{0}^{\mathbb{T}} (L_1)\to C_{0}^{\mathbb{T}} (L_2)$ decomposes in three main parts on its image, on the first part as a positive-weighted composition operator, on the second part the points in $L_2$ where the image of $T$ vanishes, and a third part formed by those points $s$ in $L_2$ such that the evaluation mapping $\delta_s\circ T$ is non-continuous. Among the consequences of this representation, we obtain that every linear bijection preserving orthogonality between commutative JB$^*$-triples is automatically continuous and biorthogonality preserving.
\end{abstract}

\maketitle

%%%%%%%%%%%%%%%%%%%%%%%%%%%%
%%%%% Introduction and Preliminaries %%%%%
%%%%%%%%%%%%%%%%%%%%%%%%%%%%

\section{Introduction}

Along this note, $\mathbb{F}$ will stand for $\mathbb{R}$ or $\mathbb{C}$.  A couple of functions $a,b$ from a set $ L$ to $\mathbb{K}$ are called \emph{disjoint} or \emph{orthogonal} if their pointwise product (denoted by $a b$) is zero, that is, they have disjoint cozero sets, where the cozero set of a function $a:L\to \mathbb{F}$ is defined by $\coz(a) = \{ t\in L : a(t) \neq 0\}$. Let $\mathcal{F}(L, \mathbb{F})$ denote the linear space of all functions from a set $L$ to $\mathbb{F}$. A mapping $\Delta : \mathcal{F}(L_1, \mathbb{F})\to \mathcal{F}(L_2, \mathbb{F})$ is called \emph{separating} or \emph{disjointness preserving} or \emph{orthogonality preserving} if it maps orthogonal or disjoint functions to orthogonal or disjoint functions (i.e. $a b =0 \Rightarrow \Delta(a) \Delta(b)$). If the equivalence $a b =0 \Leftrightarrow \Delta(a) \Delta(b) =0$  holds for all $a,b\in  \mathcal{F}(L_1, \mathbb{F})$ we say that $\Delta$ \emph{preserves orthogonality in both directions} or is \emph{biorthogonality preserving}. 

After the pioneering contributions by Abramovich \cite{Abramo83},  Beckenstein, Narici, and Todd \cite{BeckNariciTodd88}, Zaanen \cite{Zaanen75}, and Arendt \cite{Arendt83}, among others, it was revealed that for every orthogonality preserving bounded linear mapping $T: C(K_1) \to C(K_2),$ where $K_1$ and $K_2$ are compact Hausdorff spaces, there exists $h\in C(K_2)$ and a mapping $\varphi: K_2 \to K_1$ being
continuous on the set $\{ t\in K_2 : h(t) \neq 0\}$ satisfying 
$$T(f) (t) = h(t) f(\varphi (t)),$$ for all $f\in C(k)$, $t\in K$. By relaxing the hypothesis of continuity, Beckenstein, Narici and Todd found certain conditions on a orthogonality preserving linear mapping $T: C(K_1)\to C(K_2),$ and on the compact Hausdorff spaces $K_1$ and $K_2$, to guarantee that such a mapping $T$ is automatically norm continuous. That is, for example, the case if $K_1$ is first countable and if $T$ satisfies that for any pair of points $s_1\neq s_2$ in $K_2$ there exist $a,b$ in $C(K_1)$ such that $\overline{\coz(a)} \cap \overline{\coz(b)}=\emptyset$ and $T(a)(s_1)\neq 0$, $T(b) (s_2)\neq 0$. For any real or complex valued continuous function $a$ on a locally compact Hausdorff space $L$, the symbol $\supp(a)$ will stand for the support of $a$, that is, the closure in $L$ of the cozero set of $a$.\smallskip

The most conclusive and influencing result was established by Jarosz in  \cite{jarosz90}, where he provided a complete description of all linear orthogonality preserving maps between $C(K)$-spaces. A conclusion which was subsequently extended by Jeang and Wong to the setting of orthogonality preserving linear maps between $C_0 (L)$-spaces for locally compact Hausdorff spaces $L$ (see \cite{jeang96}). The main consequences of these description assert that every orthogonality preserving linear bijection between $C(K)$- or $C_0(L)$-spaces is automatically continuous and bi-orthogonality preserving. Font and Hern{\'a}ndez extended this conclusion for orthogonality preserving linear maps between subalgebras of $C_0 (L)$-spaces \cite{FontHernandez94}.

Commutative unital real C$^*$-algebras can be also represented as spaces of continuous functions. By virtue of the Gelfand theory for commutative real C$^*$-algebras, every commutative unital real C$^*$-algebra $A$ is C$^*$-isomorphic (and hence isometric) to a real function algebra of the form $C(K)^{\tau}=\{f\in C(K) : \tau (f) =f\}$, where $K$ is a compact Hausdorff space, $\tau$ is a conjugate linear $^*$-automorphism of period-2 on $C(K)$ expressed, via Banach-Stone theorem, in the form $\tau (f) (t) =\overline{f(\sigma(t))}$ ($ t\in K$), where $\sigma : K\to K$ is a topological involution (i.e. a period-2 homeomorphism) (compare \cite[Proposition 5.1.4]{LiBook2003}). In the setting of real C$^*$-algebras it is established in \cite{GarPe2014, Pe2017} that every orthogonality preserving linear bijection between commutative unital real C$^*$-algebras is automatically continuous.

In a general C$^*$-algebra $A$ the notion of ``orthogonality'' for any couple of elements $a,b\in A$ can be extended in two directions: zero product $a b =0$ and C$^*$-orthogonality $a b^* = b^* a =0$ (denoted by $a\perp b$). Working with the second notion, Wolff gave a complete characterisation of all bounded linear and symmetric
orthogonality preserving maps from a unital C$^*$-algebra into another C$^*$-algebra \cite{Wolf94}. Each C$^*$-algebra belongs to the strictly wider class of JB$^*$-triples as introduced in \cite{Ka83} with respect to the triple product $\{a,b,c\} = \frac12 (a b^* c + c b^* a)$. C$^*$-algebra problems can benefit from both the language and the techniques of JB$^*$-theory. That was the case of (C$^*$-)orthogonality preserving bounded (non-necessarily symmetric) linear operators between general C$^*$-algebras whose precise form was described by Burgos, Fern{\'a}ndez-Polo, Garc{\'e}s, Mart{\'i}nez and the second author of this note in \cite{BurFerGarMarPe08}. By providing a partial positive answer to a problem posed by Araujo and Jarosz in \cite{AraujoJarosz2003}, it was proved in  \cite{BurGarPe2011,BurGarPe2011triples} that every biorthogonality preserving linear surjection between two compact C$^*$-algebras or between two von Neumann algebras or between two atomic JBW$^*$-triples containing no infinite-dimensional rank-one summands is automatically continuous. 

Despite the multiple advances on orthogonality preserving linear maps between C$^*$- and JB$^*$-algebras (cf. \cite{BurFerGarPe09}), almost nothing is known about orthogonality preserving linear maps between general JB$^*$-triples, not even in the case of commutative JB$^*$-triples. This paper presents the first advances on orthogonality preserving linear maps between JB$^*$-triples. It is worth to note that JB$^*$-triples form a class of complex Banach spaces strictly wider than $C^*$-algebras whose roots are on the theory of holomorphic maps on complex Banach spaces of arbitrary dimension and the classification of bounded symmetric domains (the abstract definition will not be presented here, the reader can consult \cite{Ka83} and the subsequent references). In this note we shall focus on the particular subclass of all commutative JB$^*$-triples. The complex Banach spaces in this latter subclass admit a representation as closed subspaces of $C_0(L)$-spaces that are not, in general, closed for the usual pointwise product but admits a triple product $\{\cdot,\cdot,\cdot\}$ satisfying certain axioms (see \cite{Ka83} for more details). More concretely, let $\mathbb{T}$ denote the unit sphere of $\mathbb{C}$. By the Gelfand theory for JB$^*$-triples (see \cite[Corollary 1.11]{Ka83}), each commutative JB$^*$-triple $E$ can be (isometrically) identified via a triple isomorphism (i.e., a linear bijection preserving triple product), with the norm closed subspace of a $C_0(L)$ consisting of all $\mathbb{T}$-homogeneous (or $\mathbb{T}$-equivariant) continuous functions on a principal $\mathbb{T}$-bundle $L$, that is, a $\mathbb{T}$-symmetric (i.e., $\mathbb{T} L = L$) subset of a locally convex Hausdorff complex linear-space $L$ such that $0 \notin L$ and $L \cup \{0\}$ is compact with triple product given by $\{a,b,c\} = a \overline{b} c$ with pointwise products (see also \cite{cueto2021exploring}). In other words, 
\[E\cong C_0^\mathbb{T}(L):=\{a\in C_0(L):a(\lambda t)=\lambda a(t)\text{ for every } (\lambda,t)\in\mathbb{T}\times L\}.\]
The space $C_0^\mathbb{T}(L)$ is equipped with the supremum norm and the triple product defined above.

Our aim in this note is to provide a full description of all orthogonality preserving linear maps between commutative JB$^*$-triples. Our main result reads as follows:

\begin{thm}\label{thm: op jb-triples main thm}
Let $L_1$ and $L_2$ be two principal $\mathbb{T}$-bundles, and $T:C_0^\mathbb{T}(L_1)\rightarrow C_0^\mathbb{T}(L_2)$ a (non-necessarily continuous) linear orthogonality preserving map. Then, the following statements hold:
\begin{enumerate}[$(a)$] \item $L_2$ admits a decomposition as a disjoint union $L_2=L_2^z\cup L_2^d\cup L_2^c$ of (possibly empty) $\mathbb{T}$-symmetric subsets \begin{align*}
    L_2^z&:=\{s\in L_2 : \delta_s\circ T \equiv 0\}, \\
    L_2^d&:=\{s\in L_2 : \delta_s\circ T \text{ is discontinuous}\}, \\
    L_2^c&:=\{s\in L_2 : \delta_s\circ T \neq 0 \text{ and continuous}\},
\end{align*} where $L_2^z$ is closed and $L_2^d$ open. 
\item  There exists a multi-valued mapping $\hat{\varphi}$ from $L^{nz}_2=L_2^c\cup L_2^d$ into $\widehat{L_1}=L_1\cup\{0\}$ satisfying the following properties:\begin{enumerate}[$(b.1)$] \item $\hat\varphi (\lambda s) = \hat\varphi (s) $  for all $\lambda\in \mathbb{T}, s\in L^{nz}_2$.
\item For each $s\in L^{nz}_2$ and any $t_1,t_2\in \hat{\varphi} (s)$ we have $t_1 \in \mathbb{T} t_2$. The map $\hat{\varphi}$ can be viewed as a continuous mapping from $L_2^{nz}$ into ${\widehat{L_1}}/{\mathbb{T}}$, defined by $s\mapsto q_1(t),$ where $q_1$ is the canonical projection from $\widehat{L_1}$ onto $\widehat{L_1}/\mathbb{T}$ and $t$ is any element in $\widehat{\varphi}(s).$  
    \item For every $a\in C_0^\mathbb{T}(L_1)$,
\begin{equation*} \label{eq: zero if supp s not in supp a}
    \hat{\varphi}(s)\not\subset \supp(a)\Longrightarrow \delta_s\circ T(a)=T(a)(s)=0\quad\forall s\in L_2^c\cup L_2^d.
\end{equation*}
\end{enumerate}
\item There is a continuous selection of $\hat{\varphi}$ in $L_2^c$, $\varphi:L_2^c\rightarrow \widehat{L_1}$, and a non-vanishing bounded continuous function $r:L_2^c\rightarrow\mathbb{R}^+$ satisfying:  $\varphi (\lambda s) = \lambda \varphi (s)$ and $r(\lambda s) = r(s),$  for all $\lambda\in \mathbb{T}, s\in L^{c}_2$ and 
\begin{equation}\label{eq: op representation}
    T(a)(s)=r(s)\ a(\varphi(s)),\quad \hbox{for all } s\in L_2^c \ \hbox{ and every $a\in C_0^\mathbb{T}(L_1)$.}
\end{equation}
\item The set $F=\hat{\varphi}(L_2^d)$ is a disjoint finite union of the form $\mathbb{T}t_1\cup\cdots\cup\mathbb{T}t_m$ with each $t_j\in\widehat{L_1}$.
\end{enumerate}
\end{thm}

It is worth to note that in Theorem \ref{thm: op jb-triples main thm} we understand that two elements $a,b$ in a commutative JB$^*$-triple $C_0^{\mathbb{T}}(L)$ are orthogonal if they have zero pointwise product, even if the function defined by their pointwise product is not, in general, an element in $C_0^{\mathbb{T}}(L)$. This is not the formal definition of orthogonality in JB$^*$-triples. Elements $a,b$ in a JB$^*$-triple $E$ are called orthogonal if $\{a,b,x\}=0$ for all $x\in E$ (cf. \cite{BurFerGarMarPe08}). However, it can be easily checked, from the expression of the triple product on $C_0^{\mathbb{T}}(L)$ in the form $\{a,b,c\} = a \overline{b} c,$ that $a$ and $b$ are orthogonal in $C_0^{\mathbb{T}}(L)$ if and only if they have zero pointwise product, or equivalently, disjoint cozero sets.

It is known that, for each locally compact space $\tilde{L}$, the Banach space $C_0(\tilde{L})$ is isometrically isomorphic to a $C^{\mathbb{T}}_0 (L)$-space for an appropriate principal $\mathbb{T}$-bundles $L$ (cf. \cite[Proposition 10]{Ol74} or \cite[Lemma 3.1]{FriedmanRusso82TAMS}). However, there exist principal $\mathbb{T}$-bundles $L$ for which the space $C^{\mathbb{T}}_0 (L)$ is not isometrically isomorphic to a $C_0(L)$-space (cf. \cite[Corollary 1.13 and subsequent comments]{Ka83}). So, there exist abelian JB$^*$-triples which are not isometrically isomorphic to commutative C$^*$-algebras. Therefore, the previous result is a strictly sharper version of the result by Jeang and Wong for linear orthogonality preservers on $C_0(L)$-spaces (cf. \cite{jeang96}).

The main theorem will be proved along section \ref{sec:proofs} via a non-trivial adaptation of the available techniques employed by Abramovich \cite{Abramo83}, Jarosz \cite{jarosz90} and Jeang and Wong \cite{jeang96}. The arguments in the setting of $C_0^{\mathbb{T}}(L)$ will require a more sophisticated interplay between algebra, topology and analysis. Among the conclusions of our main result we shall establish in Theorem \ref{thm: op jb-triples bijective case} that every orthogonality preserving linear bijection between commutative JB$^*$-triples (i.e. $C_0^{\mathbb{T}}(L)$) is bi-orthogonality preserving and automatically continuous. 

We cannot conclude this section without commenting that $C_0^{\mathbb{T}}(L)$-spaces also appear in other research lines in functional analysis. A Banach space is called a \emph{Lindenstrauss space} if its dual is isometric to $L_1(\Omega,\Sigma,\mu)$ for some measure space $(\Omega,\Sigma,\mu)$ (cf. \cite{LindeWul1969,HirsLaz73,Ol74, FriedmanRusso82TAMS}). An isometric characterization of Lindenstrauss spaces was established by Olsen in \cite{Ol74}. A locally compact Hausdorff space $L$ is called a \emph{locally compact $\mathbb{T}_{\sigma}$-space} if there exists a continuous mapping $\mathbb{T} \times L \to  L,$ $(\lambda, t) \mapsto \lambda t$, satisfying $\lambda (\mu t) = (\lambda \mu) t$ and $1 t = t$,  for all $\lambda, \mu \in \mathbb{T}$, $t\in L$. Every principal $\mathbb{T}$-bundle is a locally compact $\mathbb{T}_{\sigma}$-space. We can extend the product by elements in $\mathbb{T}$ to the one-point compactification $L\cup \{\omega\}$ of $L$ by simply setting $\lambda \omega = \omega$ ($\lambda \in \mathbb{T}$). The Banach subspace of $C_0(L)$ of all $\mathbb{T}$-homogeneous or $\mathbb{T}$-equivariant functions on $L$ is usually denoted by $C_{\sigma} (L)$, that is,
$$C_{\sigma} (L) := \{ a \in C_0(L) : a (\lambda t) = \lambda a(t) \hbox{ for every } (\lambda, t) \in \mathbb{T}\times X\}.$$ A complex $C_{\sigma}$-space is any complex Banach space which is isometric to a space $C_{\sigma} (L)$ for some $\mathbb{T}_{\sigma}$-space $L$. A complex $C_{\Sigma}$-space is a Banach space which is isometric to a $C_{\sigma} (K)$ for some compact $\mathbb{T}_{\sigma}$-space $K$ such that $\lambda t \neq t$ for all $t\in K,$ $\lambda\in \mathbb{T}\backslash \{1\}.$ Every $C_0(L)$-space is a $C_{\sigma} (L)$-space (cf. \cite[Proposition 10]{Ol74}), and every $C_{\sigma} (L)$-space is a commutative JB$^*$-triple in the sense of Kaup \cite{Ka83} (see also \cite[Section 2]{FriedmanRusso82TAMS}). It was proved by Olsen that a complex Banach space is a $C_{\sigma}$-space (respectively, a $C_{\Sigma}$-space) if and only if it is a Lindenstrauss space and the union of $\{0\}$ with the extreme points of the closed unit ball of its dual space is weak$^*$-closed (respectively, it is a Lindenstrauss space and the set of extreme points of the closed unit ball of its dual space is weak$^*$-closed) \cite[Theorems 12 and 14]{Ol74}. Friedman and Russo proved in \cite[Theorem 5]{FriedmanRusso82TAMS} that the image of each contractive projection on $C_0(L)$ is a $C_{\sigma}$-space.

\section{Proof of the main result}\label{sec:proofs}

This section is devoted to provide a proof of Theorem \ref{thm: op jb-triples main thm}. It will be obtained after a series of results in which we mix techniques of algebra, topology and analysis. 

Henceforth, given a normed space $X$, we shall denote its unit sphere by $S(X)$. The next lemma, which has been borrowed from \cite[Remark 3.4]{cueto2021exploring} states a kind of Urysohn's lemma in the context of $C_0^\mathbb{T}(L)$-spaces. 

\begin{lem}\label{lm:T-symmetric Urysohn}{\rm\cite[Remark 3.4]{cueto2021exploring}}
Suppose $L$ is a principal $\mathbb{T}$-bundle. Let $W$ be a $\mathbb{T}$-symmetric open neighbourhood of $t_0$ in $L$ which is contained in a compact $\mathbb{T}$-symmetric subset. Then, there exists a function $h\in S(C_0^\mathbb{T}(L))$ satisfying $h(t_0)=1$ and $h(t)=0$ for all $t\in X\backslash W$.
\end{lem}

%Note that this also provides an Urysohn's lemma for $C_0^{hom}(L)$-spaces just by taking the modulus.

As the previous lemma, the next $\mathbb{T}$-symmetric separation result will be frequently employed along the document.

\begin{lem} \label{lem:disjoint T-symmetric neighbourhoods} {\rm\cite[Lemma 3.5]{cabezas22}} If $\mathbb{T}t_1\neq \mathbb{T}t_2$ in a principal $\mathbb{T}$-bundle $L$, there exist open $\mathbb{T}$-symmetric subsets $V_1,V_2\subset L$ satisfying:
\begin{enumerate}[$\bullet$]
    \item $\mathbb{T}t_j\subset V_j$ for $j=1,2$;
    \item $\overline{V_j}$ is compact for $j=1,2$;
    \item $V_1\cap V_2=\overline{V_1}\cap\overline{V_2}=\emptyset$.
\end{enumerate}
\end{lem}

In \cite[Section 3]{jeang96}, Jeang and Wong adapt the arguments by Jarosz in \cite{jarosz90} to the case of orthogonality preserving operators on $C_0(L)$-spaces using the one-point compactification. We need to go one step further in order to deal with the $\mathbb{T}$-symmetric nature of our problem.

For a principal $\mathbb{T}$-bundle $L$, we shall consider the one-point compactification $\widehat{L}=L\cup\{0\}$. Now, let us take the quotient space ${\widehat{L}}/{\mathbb{T}}$, with the equivalence relation given by the group $\mathbb{T}$, i.e.
\[t_1\sim t_2\ \Longleftrightarrow \ t_2\in\mathbb{T}t_1 \text{ for $t_1,t_2\in L$ and $\mathbb{T}0=\{0\}$}.\]
When equipped with the quotient topology, ${\widehat{L}}/{\mathbb{T}}$ is a compact Hausdorff space. The quotient map $q:\widehat{L}\rightarrow {\widehat{L}}/{\mathbb{T}}$ is always continuous and, in our case, it is also open. Namely, if $U\subset \widehat{L}$ is open, $q(U)$ is open in the quotient if and only if $q^{-1}\big(q(U)\big)=\mathbb{T}U$ is open in $\widehat{L}$, which is clear.

Given a locally compact $\mathbb{T}_{\sigma}$-space $L$, we consider the commutative C$^*$-algebra of all \emph{$\mathbb{T}$-invariant} functions in $C_0(L),$ that is,  $$C_0^{hom}({L})=\{f\in C_0({L}): f(\lambda t)=f(t), \ \forall t\in {L}, \lambda\in\mathbb{T}\},$$
equipped with the pointwise product and complex conjugation and the supremum norm. When $L$ is compact the symbol $C^{hom}({L})$ will denote the unital commutative C$^*$-algebra of all $\mathbb{T}$-invariant functions in $C(L).$  If $C_0^\mathbb{T}(L)$ is a commutative JB$^*$-triple for a principal $\mathbb{T}$-bundle $L$, the pointwise product $fa$ clearly lies in $C_0^\mathbb{T}(L)$ for any $a\in C_0^\mathbb{T}(L)$ and $f\in C^{hom}(\widehat{L})$, where $\widehat{L}=L\cup\{0\}$. In other words, $C_0^\mathbb{T}(L)$ is a $C^{hom}(\widehat{L})$-module under the pointwise product. Furthermore, the mapping $(a,b)\mapsto \langle a , b\rangle (t) := a(t)\overline{b(t)}$ defines an structure of Hilbert $C^{hom}(\widehat{L})$-module on  $C_0^\mathbb{T}(L)$ in the usual sense. 

For each $f\in C^{hom}(\widehat{L})$, the universal property of the quotient topology ensures the existence of a continuous function $\hat{f}:{\widehat{L}}/{\mathbb{T}}\rightarrow\mathbb{C}$ such that $\hat{f}(t+\mathbb{T})=f(t)$ for all $t+\mathbb{T}=\mathbb{T}t\in {\widehat{L}}/{\mathbb{T}}$.
\begin{center}
\begin{tikzcd}
    \widehat{L} \arrow{r}{f} \arrow[swap]{d}{q} & \mathbb{C}  \\
     {\widehat{L}}/{\mathbb{T}} \arrow[swap,dashed]{ru}{\hat{f}}
\end{tikzcd}
\end{center}

\begin{lem} \label{lem: hat properties} The operator \  $\widehat{\cdot}:C^{hom}(\widehat{L})\rightarrow C\big({\widehat{L}}/{\mathbb{T}}\big)$ is a well-defined linear surjective and isometric $^*$-homomorphism.
\end{lem}

\begin{proof}
The fact that  $\ \widehat{\cdot}\ $ is a linear $^*$-homomorphism follows clearly from the pointwise $^*$-operations. Next, since
\[\|\hat{f}\|_\infty=\sup\left\{|\hat{f}(t+\mathbb{T})|:t+\mathbb{T}\in {\widehat{L}}/{\mathbb{T}}\right\}=\sup\left\{|f(t)|:t\in \widehat{L}\right\}=\|f\|_\infty\]
for all $f\in C^{hom}(\widehat{L})$, $\ \widehat{\cdot}\ $ is an isometry. Thus,
$B:=\widehat{\cdot}\left(C^{hom}(\widehat{L})\right)$ is a $C^*$-subalgebra of $C\big({\widehat{L}}/{\mathbb{T}}\big)$. We shall next proof that $B$ separates points and zero in ${\widehat{L}}/{\mathbb{T}}$. Given $t_1+\mathbb{T}\neq t_2+\mathbb{T}$ in ${\widehat{L}}/{\mathbb{T}}$, we have $\mathbb{T}t_1\cap\mathbb{T}t_2=\emptyset$ in $\widehat{L}$. Thus, by Lemma \ref{lem:disjoint T-symmetric neighbourhoods}, there exist disjoint open $\mathbb{T}$-symmetric neighborhoods $V_1,V_2$ of $t_1$ and $t_2$ respectively with disjoint $\mathbb{T}$-symmetric compact closures. Using Lemma \ref{lm:T-symmetric Urysohn}, we can find orthogonal functions $a_1,a_2\in C^\mathbb{T}(\widehat{L})$ satisfying $a_j(t_j)=1$ and $a_j|_{\widehat{L}\backslash V_j}\equiv 0$ for $j=1,2$. Taking (for example) $b=a_1+2a_2\in C^\mathbb{T}(\widehat{L})$, we have $b(t_j)=j$ for $j=1,2$. Hence, $|b|\in C^{hom}(\widehat{L})$ satisfies
$0\neq \hat{|b|}(t_1+\mathbb{T})=1\neq \hat{|b|}(t_2+ \mathbb{T})=2.$

Finally, we conclude by the Stone-Weierstrass theorem \cite[Corollary 8.2]{ConwayBook} that $B=C\big({\widehat{L}}/{\mathbb{T}}\big)$, so $\ \widehat{\cdot}\ $ is surjective.
\end{proof}

We are now in condition to work within the framework of Theorem \ref{thm: op jb-triples main thm}. We recall that, as in the statement of the main theorem, we set
\begin{align*}
    L_2^z&:=\{s\in L_2 : \delta_s\circ T \equiv 0\}, \\
    L_2^d&:=\{s\in L_2 : \delta_s\circ T \text{ is discontinuous}\}, \\
    L_2^c&:=\{s\in L_2 : \delta_s\circ T \neq 0 \text{ and continuous}\}.
\end{align*}
Clearly, $L_2^z$, $L_2^d$ and $L_2^c$ are mutually disjoint, $\mathbb{T}$-symmetric and $L_2=L_2^z\cup L_2^d\cup L_2^c$. Besides, $L_2^z$ is closed. Namely, if $(s_\lambda)_\lambda$ is a net in $L_2^z$ converging to $s_0\in L_2$, for each $a\in C_0^{\mathbb{T}}(L_1)$ we have $0=\delta_{s_\lambda}\circ T(a)=T(a)(s_\lambda)\to T(a)(s_0)$ by the continuity of $T(a)$ for all $a\in C_0^\mathbb{T}(L_1)$, so $\delta_{s_0}\circ T\equiv 0$. Hence, the set $L_2^{nz}=L_2\backslash L_2^z=L_2^d\cup L_2^c$ is open in $L_2$.

For each $s\in L_2$, let us define the set $\supp(\delta_s\circ T)$ of all $t\in \widehat{L_1}$ such that for every open $\mathbb{T}-symmetric$ neighborhood $U\subset \widehat{L_1}$ of $t,$ there exists a function $a\in C_0^\mathbb{T}(L_1)$ satisfying $\coz(a)\subset U$ and $\delta_s\circ T(a)=T(a)(s)\neq 0$. It is clear that $\mathbb{T}\supp(\delta_s\circ T)=\supp(\delta_s\circ T)$.

\begin{prop} \label{prop: supp non-empty} The set $\supp(\delta\circ T)$ is empty if and only if $s\in L_2^z$ {\rm(}$\delta_s\circ T\equiv 0${\rm)}.
\end{prop}
\begin{proof}
Assume that $\supp(\delta\circ T)=\emptyset$ for some $s\in L_2$. Then, for each $t\in \widehat{L_1}$, there exists an open $\mathbb{T}$-symmetric neighborhood $U_t\subset \widehat{L_1}$ such that for all $a\in C_0^\mathbb{T}(L_1)$ with $\coz(a)\subset U_t$, $\delta_s\circ T(a)=0$. Since $\widehat{L_1}$ is compact, the open cover $\widehat{L_1}=\bigcup_{t\in \widehat{L_1}}U_t$ must have a finite subcover $\{U_1,\ldots, U_n\}$. Since the canonical projection $q: \widehat{L_1}\to \widehat{L_1}/\mathbb{T}$ is open, $\{\widehat{U_j}:=q(U_j): j=1,\ldots, n\}$ is an open cover of the (compact and Hausdorff) quotient space ${\widehat{L_1}}/{\mathbb{T}}$. Therefore, there exists a partition of unity subordinated to the cover: $\{\widehat{h_j}: j=1,\ldots, n \}\subset  C\big({\widehat{L}}/{\mathbb{T}}\big)$. By Lemma \ref{lem: hat properties}, each $\widehat{h_j}$ is the image under $\ \widehat{\cdot}\ $ of a function $h_j\in C^{hom}(\widehat{L_1})$. For all $t\in \widehat{L_1}$ we have
\[\sum_{j=1}^n h_j(t)=\sum_{j=1}^n\widehat{h_j}(t+\mathbb{T})=1.\]
Furthermore, if $t\notin U_j$, then $t+\mathbb{T}\notin \widehat{U_j}$, so $h_j(t)=\widehat{h_j}(t+\mathbb{T})=0$. This shows that $h_1,\ldots, h_n$ is a partition of unity subordinated to the cover $\{U_1,\ldots, U_n\}$.

The space $C_0^\mathbb{T}(L_1)\subset C^\mathbb{T}(\widehat{L_1})$, so it is a $C^{hom}(\widehat{L_1})$-module. Finally, given any $a\in C_0^\mathbb{T}(L_1)$, we can write $a=\sum_{j=1}^n a h_j$, and $\coz(a h_j)\subset \coz(h_j)\subset U_j$. By the definition of each $U_j$, $\delta_s\circ T(ah_j)=0$, and the linearity of $\delta_s\circ T$ leads to $\delta_s\circ T(a)=0$. The arbitrariness of $a\in C_0^\mathbb{T}(L_1)$ allows us to deduce that $s\in L_2^z$.
\end{proof}

Now, let us enumerate some properties of the set defined above.

\begin{lem} \label{lem: supp properties} The following statements hold:
\begin{enumerate}[$(1)$]
    \item\label{property: supp single orbit} $\#\ {\supp(\delta_s\circ T)}/{\mathbb{T}}= \#\ q({\supp(\delta_s\circ T)})=1,$  for all $s\in L_2^{nz}$.
    \item\label{property: supp of unitary multiple} If $s_1=\lambda s_2\in L_2$ with $\lambda\in\mathbb{T}$, then $\supp(\delta_{s_1}\circ T)=\supp(\delta_{s_2}\circ T)$.
    \item \label{property: zero if supp s not in supp a} If $\supp(\delta_s\circ T)\not\subset\overline{\coz(a)}=\supp(a)$ for some $s\in L_2^{nz}$ and $a\in C_0^\mathbb{T}(L_1)$, then $\delta_s\circ T(a)=0$.
\end{enumerate}
\end{lem}

\begin{proof} ~
\begin{enumerate}[$(1)$]
    \item Suppose that $t_1,t_2\in\supp(\delta_s\circ T)$ satisfy $\mathbb{T}t_1\neq \mathbb{T}t_2$. By Lemma \ref{lem:disjoint T-symmetric neighbourhoods}, there exist disjoint open $\mathbb{T}$-symmetric neighborhoods $V_1,V_2$ of $t_1$ and $t_2$, respectively, with disjoint $\mathbb{T}$-symmetric compact closures. For $j=1,2$, since $t_j\in \supp(\delta_s\circ T)$, there exists $a_j\in C_0^\mathbb{T}(L_1)$ satisfying $\coz(a_j)\subset U_j$ and $T(a_j)(s)\neq 0$. Therefore, $s\in \coz(T(a_1))\cap \coz(T(a_2))$. However, $\coz(a_1)\cap\coz(a_2)\subset V_1\cap V_2=\emptyset$. This contradicts the fact that $T$ preserves orthogonality, so $\mathbb{T}t_1= \mathbb{T}t_2$. This shows that $\#\ {\supp(\delta_s\circ T)}/{\mathbb{T}}= \#\ q({\supp(\delta_s\circ T)})\leq 1,$ the equality follows from Proposition \ref{prop: supp non-empty}.
    \item This follows straightforwardly from the fact that the cozero set of a function in $C_0^\mathbb{T}(L_2)$ is $\mathbb{T}$-symmetric. 
    \item Take any element $t\in\supp(\delta_s\circ T)\backslash\overline{\coz(a)}$. There must exists an open $\mathbb{T}$-symmetric neighborhood $U\ni t$ such that $U\cap\coz(a)=\emptyset$. Since $t\in\supp(\delta_s\circ T)$, there exists $b\in C_0^\mathbb{T}(L_1)$ satisfying $\coz(b)\subset U$ and $\delta_s\circ T(b)\neq 0$. Since $\coz(a)\cap\coz(b)=\emptyset$ and $T$ is orthogonality preserving, we have $\coz(T(a))\cap\coz(T(b))=\emptyset$. Therefore, $T(b)(s)\neq 0$ implies $T(a)(s)=0$.
\end{enumerate}
\end{proof}

Proposition \ref{prop: supp non-empty} along with property \eqref{property: supp single orbit} in the previous lemma allow us to define a multi-valued map:
\begin{equation} \label{eq: phi multi}
    \hat{\varphi}:L_2^{nz}\rightarrow \widehat{L_1},\quad\hat{\varphi}(s)=\supp(\delta_s\circ T).
\end{equation}

Now, fix any $s\in L_2^c$ and consider the following subspaces
\begin{align*}
    J_s&=\{a\in C_0^\mathbb{T}(L_1):\hat{\varphi}(s)\not\subset\supp(a)\} \\
    &= \{a\in C_0^\mathbb{T}(L_1) : a\hbox{ vanishes on an open $\mathbb{T}$-symmetric neighborhood of }\hat{\varphi}(s)\}, \\
    K_s&=\{a\in C_0^\mathbb{T}(L_1):a(t)=0 \ \ \forall t\in \hat{\varphi}(s)\}.
\end{align*}
Take any $a\in K_s$ and any $\varepsilon>0$, consider the closed sets
\[A_\varepsilon=\{t\in L_1: |a(t)|\geq\varepsilon\},\quad B_\varepsilon=\{t\in L_1: |a(t)|\leq\varepsilon/2\}.\]
The continuity of $a$ ensures that $\hat{\varphi}(s)\subset \operatorname{int}(B_\varepsilon)$.
Now, by Urysohn's lemma and the identification of $C^{hom}(\widehat{L_1})$ and $C\big({\widehat{L_1}}/{\mathbb{T}}\big)$ given by Lemma \ref{lem: hat properties}, we can find a function $f_\varepsilon\in S\big(C^{hom}(\widehat{L_1})\big)$ satisfying $0\leq f_\varepsilon\leq 1$, $f_\varepsilon|_{A_\varepsilon}\equiv 1$ and $f_\varepsilon|_{B_\varepsilon}\equiv 0$.

The $C^{hom}(\widehat{L_1})$-module structure assures that  $f_\varepsilon a\in C_0^\mathbb{T}(L_1)$. Moreover, we have \[\operatorname{int}(B_\varepsilon)\cap\coz(f_\varepsilon a)\subset B_\varepsilon\cap\coz(f_\varepsilon a)=\emptyset,\]
so $\hat{\varphi}(s)\not\subset\overline{\coz(f_\varepsilon a)}=\supp(f_\varepsilon a)$, and we conclude that $f_\varepsilon a\in J_s$. For all $t\in L_1$, we have
\[|f_\varepsilon(t)a(t)-a(t)|=|a(t)| |f_\varepsilon(t)-1|\leq\varepsilon,\]
since $|f_\varepsilon(t)-1|\in [0,1]$ is equal to 0 unless $a(t)<\varepsilon$. Hence, $\|f_\varepsilon a-a\|\leq\varepsilon$, and we deduce that $J_s$ is dense in $K_s$. By Lemma \ref{lem: supp properties}\eqref{property: zero if supp s not in supp a}, $J_s\subset \ker(\delta_s\circ T)$, which is closed because $\delta_s\circ T$ is a continuous linear functional (recall that $s\in L_2^c$). Thus, for any $t_0\in \hat{\varphi}(s)=t_0\mathbb{T}$,  $\ker(\delta_{t_0})=K_s\subset\ker(\delta_s\circ T)$. Since the two subspaces have codimension 1, it follows $\ker(\delta_{t_0})= \ker(\delta_s\circ T)$. Therefore, there exists a non-zero scalar $r_{t_0}(s)$ (depending on both $t_0$ and $s$) satisfying $\delta_s\circ T=r_{t_0}(s)\delta_{t_0}$, i.e.,
\begin{equation} \label{eq:multi-valued representation}
    T(a)(s)=r_{t_0}(s)\ a({t_0}),\quad \forall a\in C_0^\mathbb{T}(L_1).
\end{equation}

This equation is quite similar to the one we are looking for in the main theorem (cf. equation \eqref{eq: op representation}), and it makes clear that the function $\varphi$ we are seeking must be a continuous selection of $\hat{\varphi}$.

We shall start by studying the dependence on $t_0\in \hat{\varphi}(s)$. If we take another element in $\hat{\varphi}(s)$, it will be of the form $\mu t_0$ for some $\mu\in\mathbb{T}$. Applying equation \eqref{eq:multi-valued representation}, we get
\[r_{t_0}(s)a({t_0})=T(a)(s)=r_{\mu t_0}(s)a(\mu t_0)=\mu r_{\mu t_0}(s)a( t_0)\ \ \forall a\in C_0^\mathbb{T}(L_1),\]
and taking any $a\in C_0^\mathbb{T}(L_1)$ with $a(t_0)\neq 0$ (its existence is assured by Lemma \ref{lm:T-symmetric Urysohn}) we deduce $r_{\mu t_0}(s)=\overline{\mu} r_{t_0}(s)$ for each $s\in L_2^c$ and $t_0\in \hat{\varphi}(s)=\supp(\delta_s\circ T)$. For each $s\in L_2^c$ We shall select the unique element $t\in \hat{\varphi}(s)$ satisfying $r_t(s)>0$.

\begin{lem} \label{lm:lupdth} Let $L$ be a principal $\mathbb{T}$-bundle. Each non-zero weighted evaluation functional $\alpha\delta_s$ {\rm(}with $\alpha\in\mathbb{C}\backslash\{0\}$ and $s\in L${\rm)} on $C_0^\mathbb{T}(L)$ can be uniquely written in the form $r\delta_{\hat{s}}$, where $r\in\mathbb{R}^+$ and $\hat{s}\in \mathbb{T} s \subset L$.
\end{lem}

\begin{proof}
For all $a\in C_0^\mathbb{T}(L)$, we have
$\alpha\delta_s(a)=|\alpha|\frac{\alpha}{|\alpha|}a(s)=|\alpha|a\left(\frac{\alpha}{|\alpha|}s\right)=|\alpha|\delta_{\frac{\alpha}{|\alpha|}s}(a)$, so $\alpha\delta_s=|\alpha|\delta_{\frac{\alpha}{|\alpha|}s}$.

Next, assume that $r_1\delta_{\hat{s}_1}=r_2\delta_{\hat{s}_2}$ with $r_1,r_2\in\mathbb{R}^+$ and $\hat{s}_1,\hat{s}_2\in L$. If $\hat{s}_2\notin\mathbb{T}\hat{s}_1$, we could find (by Lemma \ref{lem:disjoint T-symmetric neighbourhoods}) an open $\mathbb{T}$-symmetric neighborhood of $\hat{s}_1$, $V$, with compact $\mathbb{T}$-symmetric closure, not containing $\hat{s}_2$. Then, Lemma \ref{lm:T-symmetric Urysohn} would ensure the existence of a function $a\in C_0^\mathbb{T}(L)$ satisfying $r_1=r_1\delta_{\hat{s}_1}(a)=r_2\delta_{\hat{s}_2}(a)=0$, which is impossible. Thus, $\hat{s}_2=\mu\hat{s}_1$ for some $\mu\in \mathbb{T}$. For any $a\in C_0^\mathbb{T}(L)$, we have $r_1 a(\hat{s}_1)=r_1\delta_{\hat{s}_1}(a)=r_2\delta_{\hat{s}_2}(a)=r_2\delta_{\mu\hat{s}_1}(a)=r_2 a(\mu\hat{s}_1)=\mu r_2 a(\hat{s}_1)$, and taking (again by Lemma \ref{lm:T-symmetric Urysohn}) $a$ such that $a(\hat{s}_1)\neq 0$, we get $r_1=\mu r_2$. But since both $r_1$ and $r_2$ are positive, $\mu$ must be 1 and the uniqueness is proved.
\end{proof}

With the previous lemma in mind, we can write the functional $r_{t_0}(s)\delta_{t_0}$ in equation \eqref{eq:multi-valued representation} uniquely in the form $r(s)\delta_{\varphi(s)}$ with $r:L_2^c\rightarrow \mathbb{R}^+$ and $\varphi:L_2^c\rightarrow \widehat{L_1}$, where $\varphi$ is a selection of $\hat{\varphi}$. Thus,
\begin{equation}\label{eq new  3.2} T(a)(s)=r(s)a(\varphi(s)),\quad \forall a\in C_0^\mathbb{T}(L_1), \ s\in L_2^c. 
\end{equation}
We have arrived at equation \eqref{eq: op representation}, but some properties of the mappings $r$ and $\varphi$ remain to be shown.

\begin{lem}\label{lem: phi r unitary multiple}
For each $s\in L_2^c$ and $\mu\in\mathbb{T}$, $\varphi(\mu s)=\mu\varphi(s)$ and $r(\mu s)=r(s)$. In particular, $\varphi(L_2^c)=\hat{\varphi}(L_2^c)$.
\end{lem}

\begin{proof}
Applying equation \eqref{eq new  3.2} with $\mu s$ and $\mu\in\mathbb{T}$, we get
\[r(\mu s) a(\varphi(\mu s))=T(a)(\mu s)= \mu T(a)(s)=\mu r(s) a(\varphi(s)).\] By Lemma \ref{lem: supp properties}\eqref{property: supp of unitary multiple}, we have $\varphi(\mu s)\in \supp(\delta_{\mu s}\circ T)=\supp(\delta_{s}\circ T)=\mathbb{T}\varphi(s)$, so $\varphi(\mu s)=\gamma \varphi(s)$ for some $\gamma\in\mathbb{T}$. Then, \[\mu r(s) a(\varphi(s))=r(\mu s) a(\varphi(\mu s))=r(\mu s) a(\gamma\varphi(s))=r(\mu s)\gamma a(\varphi(s)), \hbox{ for all } a\in C_0^\mathbb{T}(L_1),\] and by taking moduli at both sides we get $r(\mu s)=r(s)>0$ (because we may assume $a(\varphi(s))\neq 0$ by Lemma \ref{lm:T-symmetric Urysohn}). Consequently,  $\mu=\gamma$ and we conclude that $\varphi$ is $\mathbb{T}$-homogeneous.
\end{proof}

The topological properties of $\varphi$ and $r$ are quite more complicated to obtain. Thanks to Lemma \ref{lem: supp properties}\eqref{property: supp single orbit}, the assignment $s\mapsto t+\mathbb{T},$ where $t$ is any element in $\hat{\varphi}(s)$, gives a well-defined mapping from $L_2^{nz}$ to ${\widehat{L_1}}/{\mathbb{T}}$. By a little abuse of notation, we shall denote this mapping by the same symbol $\hat{\varphi}$ and we shall simply note the domain and codomain spaces. 

\begin{lem} \label{lem: phi multi-valued continuous} The mapping $\hat{\varphi}: L_2^{nz} \to {\widehat{L_1}}/{\mathbb{T}}$ is continuous.
\end{lem}

\begin{proof}
If $\hat{\varphi}$ were not continuous, there would exist a convergent net $(s_\alpha)_\alpha\to s$ in $L_2^{nz}$ such that $t_\alpha+\mathbb{T}=\hat{\varphi}(s_\alpha)\not\to\hat{\varphi}(s)=t+ \mathbb{T}$. Since ${\widehat{L_1}}/{\mathbb{T}}$ is compact, we may assume that $t_\alpha+\mathbb{T}\to t'+\mathbb{T}\neq t+\mathbb{T}$. By Lemma \ref{lem:disjoint T-symmetric neighbourhoods}, there exist open disjoint neighborhoods $U, U'\subset {\widehat{L_1}}/{\mathbb{T}}$ of $t+\mathbb{T}$ and $t'+\mathbb{T}$ respectively. Let $V=q_1^{-1}(U)$ and $V'=q_1^{-1}(U')$, where $q_1$ is the quotient mapping from $\widehat{L_1}$ onto $\widehat{L_1}/\mathbb{T}$. Then $V$ and $V'$ are open disjoint $\mathbb{T}$-symmetric neighborhoods of $t$ and $t'$ respectively. 

Since $t\in\supp(\delta_s\circ T)=\hat{\varphi}(s)$, there exists $a\in C_0^\mathbb{T}(L_1)$ such that $\coz(a)\subset V$ and $T(a)(s)\neq 0$. Keeping in mind that $s_\alpha\to s$, $T(a)$ is continuous and $t_\alpha+\mathbb{T}\to t'+\mathbb{T}$, we can find $\alpha$ large enough to satisfy $T(a)(s_\alpha)\neq 0$ and $t_\alpha+\mathbb{T}\subset U'$ (which is equivalent to $t_\alpha\in V'$).

Since $t_\alpha\in\supp(\delta_{s_\alpha}\circ T)=\hat{\varphi}(s_\alpha)$, there exists $b\in C_0^\mathbb{T}(L_1)$ satisfying $\coz(b)\subset V'$ and $T(b)(s_\alpha)\neq 0$. Finally, from the above arguments we have $\coz(a)\cap\coz(b)\subset V\cap V'=\emptyset$ and $T(a)(s_\alpha)T(b)(s_\alpha)\neq 0$, which contradicts the hypothesis that $T$ preserves orthogonality. Hence, $\hat{\varphi}$ is continuous.
\end{proof}

\begin{lem}\label{lem: r locally bounded}
The mapping $r:L_2^c\rightarrow\mathbb{R}^+$ is locally bounded at every point of $L_2^c$.
\end{lem}

\begin{proof}
Fix $s_0\in L_2^c$. The open neighborhoods of $\hat{\varphi}(s_0)=\varphi(s_0)+\mathbb{T}\in {\widehat{L_1}}/{\mathbb{T}}$ are of the form $V+\mathbb{T}$ with $\varphi(s_0)\in V$ open.
Since $\hat{\varphi}: L_2^{nz} \to {\widehat{L_1}}/{\mathbb{T}}$ is continuous, for each open and $\mathbb{T}$-symmetric neighborhood $U$ of $\varphi(s_0)$ in $\widehat{L_1}$, there exists an open $\mathbb{T}$-symmetric neighborhood $W$ of $s_0$ in $L_2^c$ such that $\hat{\varphi}(W)\subset U+\mathbb{T}$, so $\hat{\varphi}(s) \in \hat{\varphi}(W)\subset U+\mathbb{T}$, and hence $\varphi(s)\in \hat{\varphi}(s) = \supp(\delta_s\circ T)\subset  \mathbb{T} U = U,$ for all $s\in W$.

Fix now $a_0\in C_0^\mathbb{T}(L_1)$ with $a_0(\varphi(s_0))=1$, it exists by Lemma \ref{lm:T-symmetric Urysohn}. Given $\varepsilon>0$, there exits an open $\mathbb{T}$-symmetric neighborhood $U_\varepsilon$ of $\varphi(s_0)$ satisfying \[|a_0(t)-1|=|a_0(t)-a_0(\varphi(s_0))|<\varepsilon,\quad\hbox{for all } t\in U_\varepsilon.\] By the conclusion in the first paragraph, we can find $W_\varepsilon$ open and $\mathbb{T}$-symmetric such that $s_0\in W_\varepsilon$ and $\varphi(s)\in U_\varepsilon$ for all $s\in W_\varepsilon$, so
\[|a_0(\varphi(s))-1|<\varepsilon,\quad\hbox{for all } s\in W_\varepsilon.\] We have shown that the mapping $s\mapsto a_0(\varphi(s))$ is locally bounded at $s_0$. In particular, for $\varepsilon=1/2$, we have \[1/2<|a_0(\varphi(s))|<3/2,\quad\text{for all $s\in W_{1/2}$.}\]

Having in mind the expression of $T$ as a composition operator on $L_2^c$ given by equation \eqref{eq new  3.2}, we have $T(a_0)(s)=r(s)a_0(\varphi(s))$ for all $s\in W_{1/2}$, so
\[r(s)|a_0(\varphi(s))|=|T(a_0)(s)|\leq \|T(a_0)\|_\infty.\]
We conclude that $r(s)\leq 2\|T(a_0)\|$ for all $s\in W_{1/2}$.
\end{proof}

\begin{prop} \label{prop: continuity of varphi}
The mapping $\varphi:L_2^c\rightarrow\widehat{L_1}$ is continuous.
\end{prop}

\begin{proof}
Fix $s_0\in L_2^c$, and let $(s_\alpha)_\alpha$ be a net in $L_2^c$ converging to $s_0$. The net $(\varphi(s_\alpha))_\alpha$ lies in the compact space $\widehat{L_1}$, so it must have a convergent subnet. We may assume (up to a renaming of the index set) that $(\varphi(s_\alpha))_\alpha$ is convergent to $t_0\in\widehat{L_1}$. 
By the continuity of the mapping  $\hat{\varphi}: L_2^{nz} \to {\widehat{L_1}}/{\mathbb{T}}$ (cf. Lemma \ref{lem: phi multi-valued continuous}), $({\varphi} (s_\alpha)+\mathbb{T})_\alpha= (\hat{\varphi} (s_\alpha))_\alpha \to \hat{\varphi}(s_0) = \varphi(s_0) +\mathbb{T}$ in ${\widehat{L_1}}/{\mathbb{T}}$. On the other hand, $({\varphi} (s_\alpha)+\mathbb{T})_\alpha \to t_0 +\mathbb{T}$. Therefore,  $\mathbb{T} t_0 = \hat{\varphi}(s_0)$, so $t_0=\gamma\varphi(s_0)$ for some $\gamma\in\mathbb{T}$. 

By Lemma \ref{lem: r locally bounded}, there exist an open $\mathbb{T}$-symmetric neighborhood $W\in s_0$ and a constant $M>0$ satisfying $0<r(s)\leq M,$ for all $s\in W$. Therefore, $r(s_\alpha)\subset [0,M]$ for $\alpha$ large enough. By the compactness of $[0,M]$, up to taking a new subnet, we may assume that $r(s_\alpha)$ converges to some $r_0\in [0,M]$.

Fix $a_0\in C_0^\mathbb{T}(L_1)$ such that $T(a_0)(s_0)\neq 0$, this is possible because $\supp(\delta{s_0}\circ T)=\hat{\varphi}(s_0)\neq\emptyset$. Let us observe that, by \eqref{eq new  3.2}, $a_0(\varphi(s_0))\neq 0$. By a new application of \eqref{eq new  3.2}, we have 
\[r(s_\alpha)a_0(\varphi(s_\alpha))=T(a_0)(s_\alpha),\]
where the left hand side term converges to $r_0 a_0(\gamma\varphi(s_0))$ and the right hand side term converges to $T(a_0)(s_0)=r(s_0)a_0(\varphi(s_0))\neq 0$. We get $r_0\gamma a_0(\varphi(s_0))=r(s_0)a_0(\varphi(s_0))$, and dividing by $a_0(\varphi(s_0))$ we conclude that $r_0=r(s_0)$ and $\gamma=1$. Therefore, $\varphi(s_\alpha)\to\varphi(s_0)$.
\end{proof}

\begin{cor}\label{cor: r continuous}
The mapping $r:L_2^c\rightarrow \mathbb{R}^+$ is continuous.
\end{cor}

\begin{proof}
Fix $s_0\in L_2^c$, since $\varphi(s_0)\in \hat{\varphi}(s_0)= \supp(\delta_{s_0}\circ T)$, given any open $\mathbb{T}$-symmetric neighborhood of $\varphi(s_0)$, $U$, we can find $a_0\in C_0^\mathbb{T}(L_1)$ such that $\coz(a_0)\subset U$ and $T(a_0)(s_0)\neq 0$. By the continuity of $T(a_0)$, there exists an open $\mathbb{T}$-symmetric neighborhood $W$ of $s_0$ in $L_2^c$ satisfying $T(a_0)(s)\neq 0$ for each $s\in W$, and thus, by equation \eqref{eq new  3.2}, $a_0(\varphi(s))\neq 0$ for each $s\in W$. Finally, observing that $r(s)=T(a_0)(s)/a_0(\varphi(s))$ for all $s\in W$, we deduce that $r$ is continuous at $s_0$.
\end{proof}

We now have full control of $T(a),$ with $a$ arbitrary in $C_0^\mathbb{T}(L_1),$ over the elements of $L_2^c$. Our next goal is a measure of the size of the image of the multi-valued mapping $\hat{\varphi}$ over the ``discontinuous part'' $L_2^d$.

\begin{prop} \label{prop: cardinal varphi L2d}
There exist $t_1, \ldots, t_m$ in $\widehat{L_1}$ such that the image of the set $L_2^d$ under the multi-valued mapping $\hat{\varphi}$ coincides with the disjoint union $\mathbb{T} t_1 \cup \cdots \cup \mathbb{T} t_m.$ Equivalently, the mapping $\hat{\varphi}: L_2^d\to \widehat{L_1}/\mathbb{T}$ has finite image.
\end{prop}

\begin{proof}
Assume on the contrary that there exits a sequence $\{s_n \}\subset L_2^d$ satisfying that the multi-valued mapping  $\hat{\varphi}$ satisfies  $\hat{\varphi}(s_n)\neq \hat{\varphi}(s_m)$ for $n\neq m$. Pick $t_n\in \hat{\varphi}(s_n)$, we have $t_n\mathbb{T}\cap t_m\mathbb{T}=\emptyset$ for all $n\neq m$. Since the functional $\delta_{s_n}\circ T$ is discontinuous, for all $M>0$ there exits $a_M^n\in C_0^\mathbb{T}(L_1)$ such that $\|a_M^n\|\leq 1$ and $T(a_M^n)(s_n)>M$.

For each $n\in\mathbb{N}$ and any open $\mathbb{T}$-symmetric neighborhood $U_n$ of $t_n$, we can find an open $\mathbb{T}$-symmetric set $V_n$ with compact $\mathbb{T}$-symmetric closure satisfying $t_n\in V_n\subset \overline{V_n}\subset U_n$ (see \cite[Remark 3.4]{cabezas22}). Take any function $f_n\in S\big(C_0^{hom}(L_1)\big)$ such that $0\leq f_n\leq 1$, $f_n|_{\overline{V_n}}\equiv 1$ and $f_n|_{L_1\backslash U_n}\equiv 0$. Let $\mathbbm{1}\in C^{hom}(\widehat{L_1})$ denote the mapping $t\mapsto 1$. We have $a_M^n=a_M^n f_n+a_M^n(\mathbbm{1}-f_n)$, with both $a_M^n f_n$ and $a_M^n(\mathbbm{1}-f_n)$ in $C_0^\mathbb{T}(L_1)$. We have $\coz(a_M^n f_n)\subset\coz(f_n)\subset U_n$
and \[\coz(a_M^n (\mathbbm{1}-f_n))\subset\coz((\mathbbm{1}-f_n))\subset L_1\backslash \overline{V_n}\subset L_1\backslash V_n.\]
Thus, $\overline{\coz(a_M^n (\mathbbm{1}-f_n))}\subset L_1\backslash V_n,$ and we get $$\supp(a_M^n (\mathbbm{1}-f_n))\cap \supp(\delta_{s_n}\circ T) = \supp(a_M^n (\mathbbm{1}-f_n))\cap \mathbb{T} t_n \subset\supp(a_M^n (\mathbbm{1}-f_n))\cap V_n=\emptyset.$$ Lemma \ref{lem: supp properties}\eqref{property: zero if supp s not in supp a} ensures that $T(a_M^n (\mathbbm{1}-f_n))(s_n)=0$, so \begin{equation} \label{eq: T(aMn)(sn) greater than M}
    M< T(a_M^n)(s_n)=T(a_M^n f_n)(s_n)
\end{equation} by linearity. Besides, $\|a_M^n f_n\|\leq \|a_M^n\|\leq 1$.

By hypothesis, we can find a sequence $(U_n)_n$ of open $\mathbb{T}$-symmetric pairwise disjoint with $t_n\in U_n$ (apply Lemma \ref{lem:disjoint T-symmetric neighbourhoods} inductively). By the above arguments, we can also find open $\mathbb{T}$-symmetric sets $V_n$ ($n\in\mathbb{N}$) with compact $\mathbb{T}$-symmetric closures satisfying $t_n\in V_n\subset\overline{V_n}\subset U_n.$

Applying equation \eqref{eq: T(aMn)(sn) greater than M} with $M=n^4$, we obtain a sequence $(a_n)_n$ in $C_0^\mathbb{T}(L_1)$ such that $\|a_n\|\leq 1$, $\coz(a_n)\subset U_n$ and $T(a_n)(s_n)>n^4$ for all $n\in\mathbb{N}$. The $a_n$'s are mutually orthogonal and lie in the closed unit ball. Since $C_0^\mathbb{T}(L_1)$ is complete the element
\[a_0=\sum_{n=1}^\infty\frac{1}{n^2}a_n\]
belongs to $C_0^\mathbb{T}(L_1)$, as the limit of an absolutely convergent series. However, for each $n\in\mathbb{N}$ we have
\begin{align*}
    T(a_0)(s_n)=\sum_{k=1}^\infty\frac{1}{k^2}T(a_k)(s_n)=\frac{1}{n^2}T(a_n)(s_n)>n^2.
\end{align*}
This is impossible, since $|T(a_0)(s_n)|\leq \|T(a_0)\|<\infty$ for each $n\in\mathbb{N}$. We therefore conclude that $\hat{\varphi}(L_2^d)$ is finite.
\end{proof}

The following statement can be easily deduced from the arguments above.

\begin{rem} \label{rmk: sn with distinct varphi bounded norm}
Let $\{s_n\}$ be a sequence in $L_2^{nz}$ such that $\hat{\varphi}(s_n)\neq\hat{\varphi}(s_m)$ for all $n\neq m$. Then, the set $\{\|\delta_{s_n}\circ T\|:n\in\mathbb{N}\}$ must be bounded.
\end{rem}

\begin{rem} \label{rmk: r globally bounded}
Let $s_1,s_2\in L_2^c$ with $\varphi(s_2)=\lambda\varphi(s_1)$ for some $\lambda\in\mathbb{T}$. Then, for every $a\in C_0^\mathbb{T}(L_1)$, by Lemma \ref{lem: phi r unitary multiple} and \eqref{eq new  3.2}, we have \begin{align*}\delta_{s_2}\circ T(a)&=T(a)(s_2)=r(s_2)a(\varphi(s_2))=r(\lambda s_1)a(\lambda\varphi(s_1))\\ &=r(\lambda s_1)a(\varphi(\lambda s_1))=T(a)(\lambda s_1)=\delta_{\lambda s_1}\circ T(a).\end{align*}
Thus, $\delta_{s_2}\circ T=\delta_{\lambda s_1}\circ T=\lambda \delta_{s_1}\circ T$, so $\|\delta_{s_2}\circ T\|=\|\delta_{\lambda s_1}\circ T\|$. This together with the previous remark and \eqref{eq new  3.2} show that the set
\[\{\|\delta_s\circ T\|: s\in L_2^c\}=\{r(s):s\in L_2^c\}\]
must be bounded. Therefore, the mapping $r:L_2^c\rightarrow \mathbb{R}^+$ is globally bounded.
\end{rem}

Let $R:=\sup\{r(s): s\in L_2^s\}\in\mathbb{R}^+$. The unique statement which remains to be proved in Theorem \ref{thm: op jb-triples main thm} is given by next corollary. 

\begin{cor}
The set $L_2^d$ is open.
\end{cor}

\begin{proof}
For any $s\in L_2^z\cup L_2^c$, we have $\|\delta_s\circ T\|=\begin{cases} r (s)\leq R,\quad &\text{if $s\in L_2^c$} \\
0, &\text{if $s\in L_2^z$}
\end{cases}$.

Given any $s_0\in \overline{L_2^z\cup L_2^c}$, there is a net $(s_\lambda)$ in $L_2^z\cup L_2^c$ converging to $s_0$. For every $a\in C_0^\mathbb{T}(L_1)$, $|T(a)(s_\lambda)|=|\delta_{s_\lambda} \circ T(a)|\leq \|a\|_\infty R$, and by the continuity of $T(a)$ we have $|\delta_{s_0}\circ T(a)|=|T(a)(s_0)|\leq \|a\|_\infty R$. We have shown that $\delta_{s_0}\circ T$ is continuous with $\|\delta_{s_0}\circ T\|\leq R$, so $s_0\in L_2\backslash L_2^d=L_2^z\cup L_2^c$.
Therefore, $L_2^z\cup L_2^c=L_2\backslash L_2^d$ is closed.
\end{proof}

The proof of Theorem \ref{thm: op jb-triples main thm} is now completed. Our next goal will be a consequence of our main theorem, which states a result on automatic continuity of orthogonality preserving linear bijections between commutative JB$^*$-triples. 

\begin{thm} \label{thm: op jb-triples bijective case} Every orthogonality preserving linear bijection between commutative JB$^*$-triples is bi-orthogonality preserving and automatically continuous. More concretely, let $L_1$ and $L_2$ be two principal $\mathbb{T}$-bundles, and $T:C_0^\mathbb{T}(L_1)\rightarrow C_0^\mathbb{T}(L_2)$ a linear bijective orthogonality preserving map. Then, there exist a $\mathbb{T}$-homogeneous or $\mathbb{T}$-equivariant homeomorphism $\varphi:L_2\rightarrow L_1$ and a non-vanishing bounded $\mathbb{T}$-invariant continuous function $r:L_2\rightarrow\mathbb{R}^+$ such that
\begin{equation*}\label{eq: op bijective representation}
    T(a)(s)=r(s)\ a(\varphi(s)),\quad \hbox{for all } s\in L_2 \hbox{ and every $a\in C_0^\mathbb{T}(L_1)$.}
\end{equation*}
As a consequence, $T$ is bi-orthogonality preserving, that is,  $T(a)\bot T(b)$ if and only if $a\bot b$ for every $a,b\in C_0^\mathbb{T}(L_1)$.
\end{thm}

Note that the previous theorem automatically guarantees the continuity of $T$ provided that $T$ satisfies three purely algebraic conditions: linearity, bijectivity and the preservation of orthogonal elements.

\begin{proof} Our departure point will be the the conclusion of Theorem \ref{thm: op jb-triples main thm}. We shall first show that $L_2^z\neq\emptyset$. For all $s\in L_2$ there exists $b\in C_0^\mathbb{T}(L_2)$ satisfying $b(s)=1$ (cf. Lemma \ref{lm:T-symmetric Urysohn}), since $T$ is onto, we can find $a\in C_0^\mathbb{T}(L_1)$ satisfying $T(a)=b$. Thus, $T(a)(s)=b(s)=1$, so $\delta_s\circ T\neq 0$. This implies that $L_2^c=L_2\backslash L_2^d$ is closed in $L_2$ as $L_2^d$ is open. %Consequently, $L_2^c\cup\{0\}$ is closed in $\widehat{L_2}=L_2\cup\{0\}$. Namely, if $t\in \overline{L_2^c\cup\{0\}}$ and $t\neq 0$, take any compact neighborhood $U\subset L_2$ of $t$. We have $U\cap (L_2^c\cup\{0\})\neq\emptyset$, but $0\notin U$ since $U\subset L_2$ is compact. Thus, $U\cap L_2^c\neq\emptyset$ for each compact neighborhood of $t$, so $t\in \overline{L_2^c}=L_2^c$. 

We shall next prove that $\varphi(L_2^c)\subset L_1=\widehat{L_1}\backslash\{0\}$. Assume on the contrary that $\varphi(s)=0$ for some $s\in L_2^c$, by equation \eqref{eq: op representation},
\[\delta_s\circ T(a)=r(s)a(\varphi(s))=r(s)a(0)=0\]
for every $a\in C_0^\mathbb{T}(L_1)$. Then, $\delta_s\circ T\equiv 0$, so $s\in L_2^z=\emptyset$, which is impossible.

%I

We also know from Theorem \ref{thm: op jb-triples main thm} that $\hat{\varphi}(L_2^d)=\bigcup\limits_{s\in L_2^d}\supp (\delta_s\circ T)$ coincides with a disjoint union of the from $\mathbb{T}t_1\cup\cdots\cup\mathbb{T}t_m,$ for some $t_1,\ldots,t_m\in \widehat{L_1}=L_1\cup\{0\}$. We shall show next that \begin{equation}\label{eq orbits in N are isolated} \begin{aligned}
&\hbox{each orbit of the form $\mathbb{T}t_j$ with $t_j\in L_1$ must be ``non-isolated'' in $L_1$,  in the  sense }\\ 
&\hbox{that for each open $\mathbb{T}$-symmetric open set $U$ containing $t_j$ we have $U\cap \left(L_1 \backslash \mathbb{T}t_j \right)\neq \emptyset.$}
\end{aligned}
\end{equation} Suppose, on the contrary, that $\mathbb{T}t_j$ is isolated in $L_1$. Find $s_j\in L_2^d$ satisfying $\supp(\delta_{s_j}\circ T)=\mathbb{T}t_j$. Then, each function $a\in C_0^\mathbb{T}(L_1)$ with $a(t_j)=0$ also satisfies $\supp(a)\cap\hat{\varphi}(s_j)=\overline{\coz(a)}\cap \mathbb{T}t_j=\emptyset$. Namely, we can find an open and $\mathbb{T}$-symmetric neighborhood of $t_j$, $U$, such that $\mathbb{T}t_j=U\cap L_1$, so $\coz(a)\cap U=\emptyset$. Applying Lemma \ref{lem: supp properties}\eqref{property: zero if supp s not in supp a}, we get $\delta_{s_j}\circ T(a)=T(a)(s_j)=0$. We deduce that $\ker(\delta_{t_j})\subset\ker(\delta_{s_j}\circ T)$, and hence $\delta_{s_j}\circ T=\alpha \delta_{t_j}$ for some $\alpha\in\mathbb{C}$. This contradicts the fact that $\delta_{s_j}\circ T$ is discontinuous.

% II

The next step will consist in proving that the $\mathbb{T}$-symmetric set $\hat{\varphi}(L_2)\cap  L_1=\big({\varphi}(L_2^c)\cup\hat{\varphi}(L_2^d)\big)\cap L_1={\varphi}(L_2^c)\cup \big(\hat{\varphi}(L_2^d)\cap L_1\big)$ is dense in $L_1$. Otherwise, there would exist an element $t\in L_1$ and an open $\mathbb{T}$-symmetric neighborhood $W\ni t$ satisfying $W\cap \hat{\varphi}(L_2)=\emptyset$. Take, by \cite[Remark 3.4]{cabezas22}, an open $\mathbb{T}$-symmetric set $V$ with compact $\mathbb{T}$-symmetric closure satisfying $t\in V\subset \overline{V}\subset W$. By Lemma \ref{lm:T-symmetric Urysohn}, we can find $a\in C_0^\mathbb{T}(L_1)$ satisfying $a(t_0)=1$ and $a|_{L_1\backslash V}\equiv 0$. Then, $\overline{\coz(a)}\subset \overline{V}\subset W$, so $\supp(a)\cap \hat{\varphi}(L_2)=\emptyset$. Thus, $\hat{\varphi}(s) =\supp(\delta_s\circ T)\cap \supp(a)=\emptyset$, and hence Lemma \ref{lem: supp properties}\eqref{property: zero if supp s not in supp a} assures that $\delta_s\circ T(a)=T(a)(s)=0$ for all $s\in L_2$. Therefore, $T(a)\equiv 0$, which is impossible because $T$ is injective. We have therefore proved that \begin{equation}\label{eq migae of varphi is dense in L1} L_1=\overline{\hat{\varphi}(L_2)\cap L_1}=\overline{{\varphi}(L_2^c)\cup \big(\hat{\varphi}(L_2^d)\cap L_1\big).}
\end{equation}
% III
Now, we claim that $\overline{{\varphi}(L_2^c)}=L_1$. Indeed, in order to simplify the notation, let us denote $\mathcal{N}=\hat{\varphi}(L_2^d)\cap L_1$. Given any open (we can assume $\mathbb{T}$-symmetric) set $W$, by \eqref{eq migae of varphi is dense in L1}, there exists $t_W$ belonging to the intersection $\left( \varphi(L_2^c)\cup\mathcal{N}\right) \cap W$. If $t_W\in\varphi(L_2^c)$ we are done. Otherwise, $t_W\in\mathcal{N}$ and it can be written in the form $t_W=\lambda_W t_{j_W}$ for some $j_W\in\{1,\ldots, m\}$ (cf. Theorem \ref{thm: op jb-triples main thm}$(d)$). Keeping in mind that there are finitely many orbits in $\mathcal{N}$ (see again Theorem \ref{thm: op jb-triples main thm}$(d)$), there is no loss of generality in assuming that $W\cap\mathcal{N}=\mathbb{T}t_{j_W}$, since we can always replace $W$ with a smaller open $\mathbb{T}$-symmetric set. We know that the orbit $\mathbb{T}t_{j_W}$ is non-isolated in $L_1$ (cf. \eqref{eq orbits in N are isolated}), so the density of $\varphi(L_2^c)\cup\mathcal{N}$ allows us to find $t'_W\in \varphi(L_2^c)\cup\mathcal{N}$ satisfying $t'_W\in W\backslash\mathbb{T}t_{j_W}$. This leads us to $t'_W\in \varphi(L_2^c),$ and consequently, the density of $\varphi(L_2^c)$ in $L_1$ follows from the arbitrariness of $W$.

% IV

We isolate now the following property of $T$: 
\begin{equation}\label{eq Ta annihiltes on L2c implies a is zero} \hbox{For each } a\in C_0^\mathbb{T}(L_1),\hbox{ we have }  T(a)|_{L_2^c}\equiv 0 \Rightarrow T(a) = a=0.
\end{equation} Namely, if $T(a)|_{L_2^c}\equiv 0$ for some $a\in C_0^\mathbb{T}(L_1)$, it follows from \eqref{eq: op representation} that  $r(s)a(\varphi(s))=0$ for all $s\in L_2^c$, and by the properties of the function $r$ we get $a(\varphi(s))=0$ for all $s\in L_2^c$. The continuity of $a$ combined with the fact $\overline{{\varphi}(L_2^c)}=L_1$ give  $a(L_1)=a(\overline{\varphi(L_2^c)})=\{0\}$, which proves that $a=T(a)=0$.

% V

We shall prove now that $L_2=L_2^c$, equivalently, $L_2^d=\emptyset$. Assume, on the contrary, that there exists $s_0\in L_2^d$. The element $s_0$ does not belong to the closed set $L_2^c$. Thus, there exists a non-zero function $b\in C_0^\mathbb{T}(L_2)$ satisfying $b(s_0)=1$ and $b|_{L_2^c}\equiv 0$ (cf. Lemma \ref{lm:T-symmetric Urysohn} and \cite[Remark 3.4]{cabezas22}). Since $T$ is onto, we can find $a \in C_0^\mathbb{T}(L_1)$ such that $T(a)=b$. Noticing that $T(a)|_{L_2^c} = b|_{L_2^c} \equiv 0,$ the property in \eqref{eq Ta annihiltes on L2c implies a is zero} leads to $a= T(a) = b=0,$ which is impossible. Therefore, $L_2^d=\emptyset$ and $L_2=L_2^c$.

Summarizing, there exist a $\mathbb{T}$-homogeneous or $\mathbb{T}$-equivariant continuous mapping $\varphi:L_2\rightarrow L_1$ and a non-vanishing bounded $\mathbb{T}$-invariant continuous function $r:L_2\rightarrow\mathbb{R}^+$ such that 
\[T(a)(s)=r(s)a(\varphi(s)), \hbox{ for all } a\in C_0^\mathbb{T}(L_1), s\in L_2.\] It remains to prove that $\varphi$ is a homeomorphism and $T$ preserves orthogonality in both directions. For the latter, we observe that if $T(a)\bot T(c)$ for some $a,c\in C_0^\mathbb{T}(L_1)$, we have $a(\varphi(s))c(\varphi(s))=0$ for all $s\in L_2$, since $r$ is non-vanishing. The density of $\varphi(L_2)$ in $L_1$ together with the continuity of $a$ and $c$ give $a(t) \ c(t)=0$, for each $t\in L_1$, witnessing that $T^{-1}$ is also orthogonality preserving. Thus, by the arguments above, $T^{-1}$ can be written in the form $T^{-1}(b) =r_1\cdot (b\circ \varphi_1)$, $\forall b\in C_0^\mathbb{T}(L_2)$, for some non-vanishing bounded $\mathbb{T}$-invariant continuous function $r_1:L_1\rightarrow \mathbb{R}^+$ and a  $\mathbb{T}$-homogeneous or $\mathbb{T}$-equivariant continuous mapping $\varphi_1:L_1\rightarrow L_2$.

Finally, it remains to show that $\varphi_1=\varphi^{-1}$. Indeed, for each $a\in C_0^\mathbb{T}(L_1)$ and $t\in L_1$, we have
\begin{align*}
    a(t)=T^{-1}(T(a))(t)=r_1(t) T(a)(\varphi_1(t))=r_1(t) r(\varphi_1(t)) a(\varphi(\varphi_1(t))).
\end{align*}
We claim that $t':=\varphi(\varphi_1(t))=t,$ for every $t\in L_1$. Namely, if $t'\notin \mathbb{T}t$, we can find via Lemmata \ref{lem:disjoint T-symmetric neighbourhoods} and \ref{lm:T-symmetric Urysohn} an element $a\in C_0^\mathbb{T}(L_1)$ such that $a(t)\neq 0$ and $a(t')=0$, which is impossible since $r$ and $r_1$ are non-vanishing. Hence, we may therefore assume that $t'=\lambda t$ for some $\lambda\in\mathbb{T}$. In such a case, we have $a(t)=r_1(t)r(\varphi_1(t))a(\lambda t)=\lambda r_1(t)r(\varphi_1(t))a(t)$ for each  $a\in C_0^\mathbb{T}(L_1)$, taking (again by Lemma \ref{lem:disjoint T-symmetric neighbourhoods}) any $a$ such that $a(t)\neq 0$ and keeping in mind that $r$ and $r_1$ are positive-valued, we must have $\lambda\in\mathbb{R}^+$, so $\lambda=1$. Similarly, it can be proved that $\varphi_1(\varphi(s))=s,$ for all $s\in L_2$.
\end{proof}

\smallskip

\textbf{Acknowledgements}

Both authors supported by Junta de Andaluc\'{\i}a grants FQM375 and PY20$\underline{\ }$00255. Second author supported by MCIN/AEI/10.13039/501100011033/FEDER ``Una manera de hacer Europa'' project no. PGC2018-093332-B-I00, Junta de Andaluc\'{\i}a grants FQM375 and PY20$\underline{\ }$00255, and by the IMAG--Mar{\'i}a de Maeztu grant CEX2020-001105-M/AEI/10.13039/501100011033.

%%%%%%%%%%%%%%
%%% Bibliography %%%
%%%%%%%%%%%%%%


\begin{thebibliography}{99}

\bibitem{Abramo83} Y.A. Abramovich, Multiplicative representation of disjointness preserving operators, \emph{Nederl. Akad. Wetensch. Indag. Math.} \textbf{45} (1983), no. 3, 265--279.

\bibitem{AraujoJarosz2003} J. Araujo, K. Jarosz, Biseparating maps between operator algebras, \emph{J. Math. Anal. Appl.} \textbf{282} (1) (2003), 48--55. 

\bibitem{Arendt83} W. Arendt, Spectral properties of Lamperti operators, \emph{Indiana Univ. Math. J.} \textbf{32} (2) (1983), 199--215.

\bibitem{BeckNariciTodd88} E. Beckenstein, L. Narici, A.R. Todd, Automatic continuity of linear maps on spaces of continuous functions, \emph{Manuscripta Math.} \textbf{62} (1988), no. 3, 257--275. 

\bibitem{BurFerGarMarPe08} M. Burgos, F.J. Fern{\' a}ndez-Polo,
J.J. Garc{\'e}s, J. Mart{\'i}nez Moreno, A.M. Peralta,
Orthogonality preservers in C*-algebras, JB$^*$-algebras and
JB$^*$-triples, \emph{J. Math. Anal. Appl.} \textbf{348} (2008), 220--233.

\bibitem{BurFerGarPe09} M. Burgos, F.J. Fern{\'a}ndez-Polo, J.J. Garc{\'e}s, A.M. Peralta, Orthogonality preservers revisited, \emph{Asian-Eur. J. Math.} \textbf{2} (2009), 387--405. 

\bibitem{BurGarPe2011} M. Burgos, J.J. Garc{\'e}s, A.M. Peralta, Automatic continuity of biorthogonality preservers between compact C$^*$-algebras and von Neumann algebras, \emph{J. Math. Anal. Appl.} \textbf{376} (2011), no. 1, 221--230.


\bibitem{BurGarPe2011triples} M. Burgos, J.J. Garc{\'e}s, A.M. Peralta, Automatic continuity of biorthogonality preservers between weakly compact JB$^*$-triples and atomic JBW$^*$-triples, \emph{Studia Math.} \textbf{204} (2011), no. 2, 97--121.


\bibitem{cabezas22} D.Cabezas, M. Cueto-Avellaneda, D. Hirota, T. Miura, A.M. Peralta, Every commutative JB$^*$-triple satisfies the complex Mazur--Ulam property, \emph{arXiv preprint arXiv:2201.06307} (2022).

\bibitem{ConwayBook} J.B. Conway, \emph{A course in functional analysis}. Second edition. Graduate Texts in Mathematics, 96. Springer-Verlag, New York, 1990.  

\bibitem{cueto2021exploring} M. Cueto-Avellaneda, D. Hirota, T. Miura, A.M. Peralta, Exploring new solutions to Tingley's Problem for function algebras, to appear in \emph{Quaest. Math.} https://doi.org/10.2989/16073606.2022.2072787.

\bibitem{FontHernandez94} J.J. Font, S. Hern{\'a}ndez, On separating maps between locally compact spaces, \emph{Arch. Math. (Basel)} \textbf{63} (1994), 158--165.

%\bibitem{FriedmanRusoo82} Y. Friedman, B. Russo, Function representation of commutative operator triple systems, \emph{J. London Math. Soc. (2)} \textbf{27}, no. 3 (1983), 513--524.

\bibitem{FriedmanRusso82TAMS} Y. Friedman, B. Russo, Contractive projections on $C_0(K)$, \emph{Trans. Amer. Math. Soc.} \textbf{273} (1982), no. 1, 57--73.

\bibitem{GarPe2014} J.J. Garc{\'e}s, A.M. Peralta, Orthogonal forms and orthogonality preservers on real function algebras, \emph{Linear Multilinear Algebra} \textbf{62} (2014), no. 3, 275--296.

\bibitem{HirsLaz73} B. Hirsberg, A.J. Lazar, Complex Lindenstrauss spaces with extreme points, \emph{Trans. Amer. Math. Soc.} \textbf{186} (1973), 141--150. 

\bibitem{jarosz90} K. Jarosz, Automatic continuity of separating linear isomorphisms, \emph{Canad. Math. Bull.} \textbf{33}, no. 2 (1990), 139--144.

\bibitem{jeang96} J.-S. Jeang, N.-C. Wong, Weighted composition operators of $C_0(X)$'s, \emph{J. Math. Anal. Appl.} \textbf{201} (1996), no. 3, 981--993.

\bibitem{Ka83} W. Kaup, A Riemann Mapping Theorem for bounded symmentric domains in complex Banach spaces, \emph{Math. Z.} \textbf{183} (1983), 503--529.

\bibitem{LiBook2003} B.R. Li, \textit{Real operator algebras}, World Scientific Publishing Co. (Singapore), Inc., River Edge, NJ, 2003.

\bibitem{LindeWul1969} J. Lindenstrauss, D.E. Wulbert, On the classification of the Banach spaces whose duals are $L_1$ spaces, \emph{J. Funct. Anal.} \textbf{4} (1969), 332--349. 

\bibitem{Ol74} G.H. Olsen, On the classification of complex Lindenstrauss spaces, \emph{Math. Scand.} \textbf{35} (1974), 237--258.

\bibitem{Pe2017} A.M. Peralta, Orthogonal forms and orthogonality preservers on real function algebras revisited, \emph{Linear Multilinear Algebra} \textbf{65} (2017), no. 2, 361--374. 

\bibitem{Wolf94} M. Wolff, Disjointness preserving operators in C$^*$-algebras, \emph{Arch. Math.} \textbf{62} (1994), 248--253.


\bibitem{Zaanen75} A.C. Zaanen, Examples of orthomorphisms, \emph{J. Approximation Theory} \textbf{13} (1975), 192--204.

\end{thebibliography}
\end{document}